# On the images of modular and geometric three-dimensional Galois representations

Luis Dieulefait*

Nria Vila†
Dept. d'lgebra i Geometria, Universitat de Barcelona;
Gran Via de les Corts Catalanes 585; 08007 - Barcelona; Spain.
e-mail: luisd@mat.ub.es; vila@mat.ub.es.

**Abstract**

We study compatible families of three-dimensional Galois representations constructed in the tale cohomology of a smooth projective variety. We prove a theorem asserting that the images will be generically large if certain easy to check conditions are satisfied. We only consider representations with coefficients in an imaginary quadratic field. For primes inert in this field, the residual representations (when irreducible) are unitary. We apply our result to an example constructed by van Geemen and Top (see [vG-T1]), obtaining a family of special linear groups and one of special unitary groups as Galois groups over $\mathbb{Q}$.
We also consider the case of a cohomological modular form for a congruence subgroup of $SL(3, \mathbb{Z})$. Assuming Clozel's conjecture stating that a geometric family of three-dimensional Galois representations can be attached to it, we give conditions on the modular form guaranteeing the validity of the result on the largeness of the images. We apply this result to several examples.

(Running head: Images of three-dimensional Galois representations)

## 1 Introduction

In this and the next section, we give an account of the theory of modular forms for congruence subgroups of $SL(3, \mathbb{Z})$, Clozel's conjecture predicting a family of Galois representations attached to Hecke eigenforms, and van Geemen and Top's construction of three-dimensional compatible families of geometric Galois representations that experimentally support this conjecture. This provides the motivation and background for the following sections, where we determine the images of these geometric representations and apply this to the realization of three-dimensional special linear and unitary groups over finite fields as Galois groups over $\mathbb{Q}$. See section 3 for more details.

*Supported by TMR - Marie Curie Fellowship ERBFMBICT983234
†Partially supported by MCYT grant BFM2000-0794-C02-01

In section 12, assuming Clozel's conjecture in a stronger form, we obtain examples of Galois representations over different fields, and an application of the techniques described in previous sections shows that for eigenforms satisfying some easy to check standard conditions the images of these Galois representations give special linear or unitary groups for all but finitely many primes.

References for the first two sections: [A-G-G] , [A-M], [vG-T2], [vG-T1] and [vG-K-T-V].

## 1.1 Modular forms for congruence subgroups of $\mathrm{SL}(3,\mathbb{Z})$

For classical modular forms, the space $S_2(\Gamma_0(N))$ can be thought of as a subspace of $H^1(\Gamma_0(N),\mathbb{C})$. This inspires an analogous definition of modular forms and Hecke operators in the three dimensional case, i.e., for congruence subgroups of $\mathrm{SL}(3,\mathbb{Z})$. For $N \geq 1$ an integer, let $\Gamma_0(N) \subseteq \mathrm{SL}(3,\mathbb{Z})$ be the subgroup whose elements $M = (m_{ij})$ satisfy:

$$m_{21} \equiv m_{31} \equiv 0 \pmod{N}.$$

**Definition 1.1** *A modular form of level $N$ is an element of $H^3(\Gamma_0(N),\mathbb{C})$.*

**Definition 1.2** *For every $\alpha \in \mathrm{GL}(3,\mathbb{Q})$, let $\Gamma_\alpha := \Gamma \cap \alpha^{-1}\Gamma\alpha$, where $\Gamma := \Gamma_0(N)$. $\Gamma_\alpha$ has finite index both in $\Gamma$ and in $\alpha^{-1}\Gamma\alpha$. We have two natural maps from $\Gamma_\alpha$ to $\Gamma$:*

$$\varphi : \gamma \to \gamma,$$
$$\psi : \gamma \to \alpha\gamma\alpha^{-1}.$$

*Consider the restriction $\varphi^* : H^3(\Gamma,\mathbb{C}) \to H^3(\Gamma_\alpha,\mathbb{C})$ and the transfer induced by $\psi$: $\psi_* : H^3(\Gamma_\alpha,\mathbb{C}) \to H^3(\Gamma,\mathbb{C})$.*
*We call Hecke operator the $\mathbb{C}$-linear operator:*

$$T_\alpha : H^3(\Gamma_0(N),\mathbb{C}) \to H^3(\Gamma_0(N),\mathbb{C})$$

*given by the composition $\psi_* \circ \varphi^*$.*

**Proposition 1.3** *Let $\mathcal{T}(N)$ be the subalgebra of $\mathrm{End}(H^3(\Gamma_0(N),\mathbb{C}))$ generated by the Hecke operators $T_\alpha$ with $\alpha \in \mathrm{GL}(3,\mathbb{Z})^+, (\det(\alpha),N) = 1$. Then $\mathcal{T}(N)$ is commutative and it is generated by the operators $T_{\alpha_p}$ and $T_{\beta_p}$ with $p \nmid N$ and:*

$$\alpha_p := \begin{pmatrix} p & 0 & 0 \\ 0 & 1 & 0 \\ 0 & 0 & 1 \end{pmatrix}; \quad \beta_p := \begin{pmatrix} p & 0 & 0 \\ 0 & p & 0 \\ 0 & 0 & 1 \end{pmatrix}.$$

$H^3(\Gamma_0(N), \mathbb{C})$ can be decomposed as direct sum of common eigenspaces for the operators in $\mathcal{T}(N)$ :

$$H^3(\Gamma_0(N), \mathbb{C}) = \bigoplus_\lambda V_\lambda,$$

with $\lambda$ homomorphism of algebras $\mathcal{T}(N) \to \mathbb{C}$ and $Tf = \lambda(T)f$, for $T \in \mathcal{T}(N), f \in V_\lambda$.

Given a character $\lambda$ of $\mathcal{T}(N)$, we define the eigenvalue $a_p$ of the Hecke eigenform $f \in V_\lambda$ :

$$a_p := \lambda(T_{\alpha_p}).$$

Then it is known that: $\lambda(T_{\beta_p}) = \bar{a}_p$ .
The field generated by the eigenvalues $a_p$ of a Hecke eigenform $f$ is a number field, either totally real or a CM-field. We will denote this field $\mathbb{Q}_f$. The eigenvalues $a_p$ are algebraic integers.

**Conjecture 1.4** *(Clozel) Let $f \in H^3(\Gamma_0(N), \mathbb{C})$ be a cuspidal eigenform for the action of the Hecke algebra $\mathcal{T}(N)$ with eigenvalues: $T_{\alpha_p}f = a_p f$ for every $p \nmid N$ (and therefore $T_{\beta_p}f = \bar{a}_p f$) generating the number field $\mathbb{Q}_f$.*
*Then, if we take for every rational prime $\ell$ a prime $\lambda$ in $\mathbb{Q}_f$ dividing $\ell$, there exists a compatible family of $\lambda$-adic Galois representations*

$$\rho_\lambda : G_\mathbb{Q} \to \mathrm{GL}(3, \mathbb{Q}_{f,\lambda})$$

*such that for every $\lambda$, $\rho_\lambda$ is unramified outside $\ell N$ and for every $p \nmid \ell N$ the characteristic polynomial of $\rho_\lambda(\mathrm{Frob}\ p)$ is:*

$$x^3 - a_p x^2 + p\bar{a}_p x - p^3.$$

See [C], [vG-T1] and [vG-K-T-V]. For the definition of cuspidality see [A-G-G].
Remark: In [A-M], a similar conjecture for mod $p$ modular forms is stated.

We will only consider the case of $\mathbb{Q}_f$ a CM-field. In the other case the family $\rho_\lambda$ becomes selfdual, and there is a finite index subgroup of the image of these representations contained in $G \cong \mathrm{PGL}(2, \mathbb{Q}_{f,\lambda})$.
Most of the known examples of eigenforms $f \in H^3(\Gamma_0(N), \mathbb{C})$ have imaginary quadratic $\mathbb{Q}_f$, such as $\mathbb{Q}_f = \mathbb{Q}(i)$ for $f \in H^3(\Gamma_0(128), \mathbb{C})$.

## 2 Experimental verification of Clozel's conjecture

A compatible family of three-dimensional $\lambda$-adic Galois representations was constructed by van Geemen and Top (see [vG-T1]) to give experimental evidence for the previous conjecture for three eigenforms of different levels, all with $\mathbb{Q}_f = \mathbb{Q}(i)$. These representations are geometric, and we summarize their construction:

Let $X$ be an algebraic variety defined over $\mathbb{Q}$, and consider the Galois $\ell$-adic representations on tale cohomology:

$$\rho_\ell : G_\mathbb{Q} \to \mathrm{GL}(H_{et}^n(X_{\overline{\mathbb{Q}}}, \mathbb{Q}_\ell)).$$

If $X$ is smooth, projective, and with good reduction at $p$, it is known that the representations $\rho_\ell$ ($\ell \neq p$), are unramified at $p$ and the characteristic polynomial of $\rho_\ell(\mathrm{Frob}\ p)$ has coefficients in $\mathbb{Z}$, is independent of $\ell$ and its roots have absolute value $p^{n/2}$.

Clozel's conjecture predicts the determinant of the representations attached to modular forms to be $p^3$, so we must take $n = 2$. This also implies $|a_p| \leq 3p$.

The examples are constructed using a subspace of $H_{et}^2(S_{\overline{\mathbb{Q}}}, \mathbb{Q}_\ell)$, where $S$ is a surface. The Galois representation decomposes:

$$H_{et}^2(S_{\overline{\mathbb{Q}}}, \mathbb{Q}_\ell) = T_\ell \oplus \mathrm{NS}(S_{\overline{\mathbb{Q}}}) \otimes_\mathbb{Z} \mathbb{Q}_\ell,$$

where $\mathrm{NS}(S_{\overline{\mathbb{Q}}})$ is the Nron-Severi group and $T_\ell$ its orthogonal complement with respect to the cup product: $H_{et}^2 \times H_{et}^2 \to H_{et}^4 \cong \mathbb{Q}_\ell$.

Eigenvalues of Frobenius on $\mathrm{NS}(S_{\overline{\mathbb{Q}}})$ are roots of unity multiplied by $p$, so we will look for our family of $\lambda$-adic representations inside $T_\ell \otimes_{\mathbb{Q}_\ell} K_\ell$ with $K = \mathbb{Q}(i)$.

If it were dim $T_\ell = 3$ we would have selfdual representations given by three dimensional $\mathbb{Q}_\ell$-vector spaces. We look for three dimensional $\lambda$-adic representations $\rho_{f,\lambda}$ associated as in the conjecture to a modular form $f$ with $\mathbb{Q}_f = \mathbb{Q}(i)$. Therefore, we want $\mathrm{trace}(\rho_{f,\lambda}(\mathrm{Frob}\ p)) \in \mathbb{Z}[i]$ and not all these traces in $\mathbb{Z}$.

Thus we suppose that the surface has an automorphism defined over $\mathbb{Q}$:

$$\phi : S \to S \quad \text{with} \quad \phi^4 = id_S.$$

The induced map: $\phi^* : H_{et}^2 \to H_{et}^2$ commutes with the Galois representations. Suppose that dim $T_\ell = 6$ and $\phi^* : T_\ell \to T_\ell$ has two three dimensional eigenspaces $W_\lambda, W_{\lambda'}$, with eigenvalues $i, -i$ respectively, such that:

$$T_\lambda := T_\ell \otimes_{\mathbb{Q}_\ell} \mathbb{Q}_\ell(i) = W_\lambda \oplus W_{\lambda'}.$$

This decomposition gives us a compatible family of $\lambda$-adic three dimensional Galois representations $\sigma'_\lambda$ on $W_\lambda$, a $\mathbb{Q}(i)_\lambda$ vector space. In general $\sigma'_\lambda$ has to be twisted by a Dirichlet character $\epsilon$ to obtain $\epsilon \otimes \sigma'_\lambda = \sigma_\lambda$ with $\det(\sigma_\lambda(\mathrm{Frob}\ p)) = p^3$.

These are the representations we will consider:

$$\sigma_{S,\lambda} : G_\mathbb{Q} \to \mathrm{GL}(W_\lambda).$$

**The Example**: The family of surfaces constructed by van Geemen and Top is the following: $S_a$ is the smooth, minimal projective model of the singular affine surface defined by:

$$t^2 = xy(x^2 - 1)(y^2 - 1)(x^2 - y^2 + axy).$$

The automorphism $\phi$ is given by $(x, y, t) \to (y, -x, t)$.

Let $N = N(a)$ be the product of the primes of bad reduction of $S_a$.

The characteristic polynomials of $\sigma_{S_a,\lambda}(\text{Frob } p)$ will be denoted

$$x^3 - a_p x^2 + p\bar{a}_p x - p^3.$$

In particular, we put $a_p := \text{trace}(\sigma_{S_a,\lambda}(\text{Frob } p)) \in \mathbb{Q}(i)$ for every $p \nmid \ell N$. These traces are determined using Lefschetz trace formula and counting points of $S_a$ over finite fields. Of course these surfaces verify $\dim T_\ell = 6$.

We will illustrate our study of three-dimensional geometric Galois representations with the examples given by the surfaces $S_a$. This is just to simplify the exposition; the results obtained can be applied to any similar example (see section 11).

In [vG-K-T-V] the following result is proved:

**Theorem 2.1** *Clozel's conjecture (conjecture 1.4) for three examples of modular forms in $H^3(\Gamma_0(N), \mathbb{C})$ with $N = 128, 160, 205$ is experimentally supported by the representations $\sigma_{S_a,\lambda}$ described above, with $a = 2, 1, 1/16$, respectively. More precisely, the eigenvalues of these eigenforms agree with the traces of the images of Frobenius for the corresponding representation, for all $p \nmid 2N$, $p \leq 173$.*

Remarks: For example, in the case $a = 2$, we have to twist by

$$\epsilon : G_\mathbb{Q} \to \text{Gal}(\mathbb{Q}(\sqrt{-2})/\mathbb{Q}) \cong \pm 1,$$

to obtain $\det(\sigma_{S_2,\lambda}) = \chi^3$, $\chi$ the $\ell$-adic cyclotomic character.

Observe that the properties possessed by a geometric representation that we will use later are preserved by this twist (description of the restriction to the inertia subgroup at $\ell$ and absolute value of roots of characteristic polynomials), so no additional difficulties are introduced.

The primes of bad reduction of $S_a$, which are the only primes other than $\ell$ where $\sigma_{S_a,\lambda}$ may ramify, are the primes dividing $2a(a^2 + 4)$. This gives only the prime 2 in the case $a = 2$. It is not known whether the conductor of the representations $\sigma_{S_a,\lambda}$ is bounded independently of $\ell$. If this is the case, we call this minimal uniform bound "conductor of the family of representations". The existence of this conductor is conjectured.

## 3 Statement of purpose

We will consider the compatible family of Galois representations $\sigma_{S_2,\lambda}$ constructed in the previous section from the surface $S_2$.

We will prove that the image is a special linear group (up to a twist by the cyclotomic character) for the primes decomposing in $\mathbb{Q}(i)$, i.e., $\ell \equiv 1 \pmod 4$, except at most for a set of Dirichlet density 0, and the image of the residual Galois representation is a unitary group for the inert primes $\ell \equiv 3 \pmod 4$, except at most for a set of Dirichlet density 0.

The exceptional density 0 sets, both for decomposing and inert primes, could be replaced by finite sets if we merely knew the existence and value of the conductor of the family $\sigma_{S_2,\lambda}$.

The method works in general for geometric three-dimensional representations over a quadratic imaginary field such as those associated to the surfaces $S_a$, whenever certain general conditions that we will describe later are satisfied by the traces of Frobenius (see section 6 and 8) and the action of the inertia subgroup at $\ell$ (see section 5).

The family of representations $\sigma_{S_a,\lambda}$ comes from the action of $G_\mathbb{Q}$ on $W_\lambda$, the 3-dimensional $i$-eigenspace (for $\phi^* : T_\ell \to T_\ell$) contained in $T_\ell \subseteq H^2_{et}(S_{a\overline{\mathbb{Q}}}, \mathbb{Q}_\ell)$. $T_\ell$ also admits a realization $T_\mathbb{Z} \subseteq H^2(S_a(\mathbb{C}), \mathbb{Z})$ which is a Hodge structure of weight 2, i.e.: $T_\mathbb{Z}$ is a $\mathbb{Z}$-module, free of rank 6 with

$$T_\mathbb{Z} \otimes \mathbb{C} = \bigoplus_{p+q=2} T^{p,q} \qquad \overline{T}^{p,q} = T^{q,p},$$

with $\dim T^{p,q} = 2$ and $\phi^*$ has eigenvalues $i, -i$ on each $T^{p,q}$. In particular $W_\lambda \otimes \mathbb{C}$ intersects non-trivially each $T^{p,q}$ (cf. [vG-T1]).

## 4  Classification of subgroups of $\mathrm{PSL}(3, \mathbb{F}_q)$

What follows is the classification, due to Mitchell ([M]), of the maximal subgroups $\mathcal{G}$, up to conjugation, of $\mathrm{PSL}(3, \mathbb{F}_{p^t})$, for $t = 1$ or $2$, as groups of collineations of the projective plane over $\mathbb{F}_{p^t}$. See [V-Y] for the relevant definitions from plane projective geometry. We include the order of the subgroup $\mathcal{G}$ in each case.

a) the stabilizer of a point: $|\mathcal{G}| = (p^t + 1)p^{3t}(p^t - 1)^2 / \gcd(3, p^t - 1)$
b) the stabilizer of a line: $|\mathcal{G}| = (p^t + 1)p^{3t}(p^t - 1)^2 / \gcd(3, p^t - 1)$
c) groups that leave invariant a triangle with coordinates in $\mathbb{F}_{p^t}$ and make all six permutations of its vertices: $|\mathcal{G}| = 6(p^t - 1)^2 / \gcd(3, p^t - 1)$
d) groups that leave invariant a triangle with coordinates in $\mathbb{F}_{p^{3t}} \smallsetminus \mathbb{F}_{p^t}$ and permute its vertices cyclically: $|\mathcal{G}| = 3(p^{2t} + p^t + 1) / \gcd(3, p^t - 1)$
e) the stabilizer of a conic: $|\mathcal{G}| = (p^t + 1)p^t(p^t - 1)$
f) $\mathrm{PSL}(3, \mathbb{F}_{p^k})$ with $k \leq t$:
$|\mathcal{G}| = (p^{2k} + p^k + 1)(p^k + 1)p^{3k}(p^k - 1)^2 / \gcd(3, p^k - 1)$
g) (only if $t = 2$) group of special unitary transformations $\mathrm{PSU}(3, \mathbb{F}_p)$: $|\mathcal{G}| = (p^2 - p + 1)(p + 1)^2 p^3 (p - 1) / \gcd(3, p + 1)$
h) Hessian groups of orders 216 (if $p^t \equiv 1 \pmod{9}$), 72 or 36 (if $p^t \equiv 4, 7 \pmod{9}$)
i) group isomorphic to $\mathrm{PSL}(2, 7)$: $|\mathcal{G}| = 168$ (if $t = 2$ or if $t = 1$, $p \equiv 0, 1, 2, 4 \pmod{7}$)
j) group isomorphic to $A_6$: $|\mathcal{G}| = 360$ (if $t = 2$ or if $t = 1$, $p \equiv 1, 4 \pmod{15}$).

We have excluded the case $p = 5$ and $t = 2$. In this case there are more possibilities.

Remarks: In cases a) and b) the subgroup is reducible. Cases c) and d) correspond to the normalizer of the group of diagonal matrices and that of the group of units of the extension $[\mathbb{F}_{p^{3t}} : \mathbb{F}_{p^t}]$ embedded in $\mathrm{SL}(3, \mathbb{F}_{p^t})$ using the regular representation respectively. In particular, inside such a group we have groups that are reducible over $\mathbb{F}_{p^t}$ (over $\mathbb{F}_{p^{3t}}$), namely, a group of diagonal matrices and a cyclic group of units embedded by the regular representation respectively. In case h), the group $\mathcal{G}$ fits in the exact sequence

$$1 \to \mathcal{F} \to \mathcal{G} \to \mathcal{U} \to 1,$$

with $\mathcal{U}$ isomorphic to a subgroup of $\mathrm{SL}(2,3)$ or the quaternionic group $\mathcal{D}_8$ and $\mathcal{F}$ of type $(3,3)$ (cf. [B]).
We will call cases h), i) and j) "exceptional cases".

## 5 Comparison of tale and Crystalline cohomologies and the action of inertia

As a consequence of the result of Fontaine-Messing in [F-M] giving (canonical, functorial) isomorphism between the tale cohomology and the crystalline cohomology for a variety $S$ having good reduction at $\ell$ and with $\dim(S) < \ell$, the following theorem, conjectured by Serre in [S2], holds:

**Theorem 5.1** *Let $V$ be a $G_\mathbb{Q}$-stable lattice in $H^m_{et}(S_{\overline{\mathbb{Q}}}, \mathbb{Q}_\ell)$ and $\tilde{V}$ the semi-simplification of $V/\ell V$ with respect to the action of $I = I_\ell$. Then, if $\dim(S) < \ell$ and $S$ have good reduction at $\ell$, the action of $I$ is given by a character $\Psi : I_t \to \mathbb{F}_q^*$ with $q = \ell^n, n = \dim V$, such that*

$$\Psi = \Psi_f^{d_0 + \ell d_1 + \ldots + \ell^{n-1} d_{n-1}}$$

*with $\Psi_f$ fundamental character of level $n$ and the exponents verify: $0 \leq d_i \leq m$.*

We apply this result to the representations $\sigma_{S_a, \lambda}$. For any $\lambda$-adic representation $\sigma_\lambda$, we denote by $\bar{\sigma}_\lambda$ the corresponding residual representation obtained by composing $\sigma_\lambda$ with the reduction modulo $\lambda$.
Since $\det(\sigma_{S_a,\lambda}) = \chi^3$ and $m = 2$, we have the following possibilities for $\bar{\sigma}_{S_a,\lambda}|_I$, where $\ell > 2$ is a prime of good reduction for $S_a$ :

$$\begin{pmatrix} 1 & * & * \\ 0 & \chi & * \\ 0 & 0 & \chi^2 \end{pmatrix}; \quad \begin{pmatrix} 1 & * & * \\ 0 & \psi_2^{2+\ell} & 0 \\ 0 & 0 & \psi_2^{1+2\ell} \end{pmatrix}; \quad \begin{pmatrix} \chi & * & * \\ 0 & \psi_2^2 & 0 \\ 0 & 0 & \psi_2^{2\ell} \end{pmatrix};$$

$$\begin{pmatrix} \chi^2 & * & * \\ 0 & \psi_2 & 0 \\ 0 & 0 & \psi_2^\ell \end{pmatrix}; \quad \begin{pmatrix} \psi_3^{1+2\ell} & 0 & 0 \\ 0 & \psi_3^{\ell+2\ell^2} & 0 \\ 0 & 0 & \psi_3^{\ell^2+2} \end{pmatrix}; \quad \begin{pmatrix} \chi & * & * \\ 0 & \chi & * \\ 0 & 0 & \chi \end{pmatrix}$$

with $\psi_i$ a fundamental character of level $i$.

We want to compare this description with the classification of the previous section. To do so, we first need to exclude the last case, the case with trivial $\mathbb{P}((\bar{\sigma}_\lambda|_I)^{ss})$.

In the examples $\bar{\sigma}_{S_a,\lambda}$ this can be done because of the Hodge structure of $T_\mathbb{Z}$. Using again the fact that our representations are crystalline ([F-M]), we are free to apply the results of Fontaine-Laffaille (cf. [F-L], theorem 5.3), so the knowledge of the Hodge filtration allows us to compute the possible values of the exponents of the fundamental characters describing the action of $I_\ell$ as in theorem 5.1. In particular the exponent 1 cannot appear three times because dimension $T^{p,q} = 2$ for $(p,q) = (0,2), (1,1)$ and $(2,0)$.

In the other cases we have $\mathbb{P}(\bar{\sigma}_{S_a,\lambda}|_I)$ of the form:

$$\begin{pmatrix} \chi^{-1} & * & * \\ 0 & 1 & * \\ 0 & 0 & \chi \end{pmatrix}; \quad \begin{pmatrix} \chi^{-1} & * & * \\ 0 & \psi_2 & 0 \\ 0 & 0 & \psi_2^\ell \end{pmatrix}; \quad \begin{pmatrix} 1 & * & * \\ 0 & \psi_2^{1-\ell} & 0 \\ 0 & 0 & \psi_2^{\ell-1} \end{pmatrix};$$

$$\begin{pmatrix} \psi_3^{\ell-\ell^2} & 0 & 0 \\ 0 & \psi_3^{\ell^2-1} & 0 \\ 0 & 0 & \psi_3^{1-\ell} \end{pmatrix}.$$

# 6 Determination of the images: totally decomposing primes

## 6.1 Reducible representations

Recall that $N$ denotes the product of the primes of bad reduction of $S_a$. Suppose that $\ell \equiv 1 \pmod 4$ is a prime such that $\bar{\sigma}_{S_a,\lambda}$ is not irreducible over $\mathbb{F}_{\ell^3}$. In this case there is a character $\mu$ unramified outside $N\ell$ with image in $\mathbb{F}_{\ell^3}^*$ such that $\mu(p)$ is a root of $x^3 - a_p x^2 + p\bar{a}_p x - p^3$, for every $p \nmid \ell N$.
Using the description of $\bar{\sigma}_{S_a,\lambda}|_I$ we know that $\mu = \chi^i \varepsilon$, with $\chi$ the mod $\ell$ cyclotomic character, $i = 0, 1, 2$ and $\varepsilon$ a character unramified outside $N$.

The problem is that we cannot bound the value of $\mathrm{cond}(\varepsilon)$, the conductor of $\varepsilon$, independently of $\ell$. We have to exclude a set of primes of small density in order to bound it.

For the family of Galois representations $\bar{\sigma}_{S_2,\lambda}$ of the example, we have $N = 2$ and $\varepsilon : (\mathbb{Z}/2^u\mathbb{Z})^* \to \mathbb{F}_{\ell^3}^*$, with $\ell \equiv 1 \pmod 4$.
We exclude the set $\mathcal{P}$ of primes $\ell \equiv 1 \pmod{64}$ in order to bound the conductor of $\varepsilon$.
What follows only holds for $\ell \equiv 1 \pmod 4$, $\ell \notin \mathcal{P}$.

With this restriction $64 \nmid \ell^3 - 1$ holds, so the image of $\varepsilon$ is a cyclic group of order at most 32. Then the conductor of $\varepsilon$ is at most 128.

Thus, $\varepsilon(p) = 1$, for every $p \equiv 1 \pmod{128}$. Moreover, $\varepsilon(p) = 1$ for every $p \equiv -1 \pmod{128}$, except if $\varepsilon$ is a character of conductor 4 or 8, in which case $\varepsilon(3) = \pm 1$.

We will use only $a_3$ and $a_{127}$ because either $\varepsilon(127) = 1$ or $\varepsilon(3) = \pm 1$. We are assuming that $\mu(p) = p^i \varepsilon(p)$, $i = 0, 1$ or $2$, is a root of the characteristic polynomials of $\bar{\sigma}_{S_2,\lambda}(\text{Frob } p)$, for every $p \nmid 2\ell$. Taking in particular $p = 127$ or $p = 3$, we know that either $\varepsilon(127) = 1$, thus giving the equation:

$$127^{3i} - a_{127}127^{2i} + 127\bar{a}_{127}127^i - 127^3 \equiv 0 \pmod{\ell},$$

$i = 0, 1$ or $2$, with $a_{127} = 161 - 16i$; or $\varepsilon(3) = \pm 1$, giving a similar equation for $p = 3$, $a_3 = 1 + 2i$, except for some $\pm$ signs (the values of the traces $a_p$ that we use are listed in [vG-K-T-V]).

From these equations we check that $\bar{\sigma}_{S_2,\lambda}$ may reduce over $\mathbb{F}_{\ell^3}$ only for the primes $\ell = 5, 13, 1901, 105649$. Looking at the characteristic polynomials of $\bar{\sigma}_{S_2,\lambda}(\text{Frob } p)$ for these $\ell$ and several primes $p$, we see that the representation $\bar{\sigma}_{S_2,5}$ seems to be reducible (over $\mathbb{F}_5$). For the other three primes it is irreducible (over $\mathbb{F}_{\ell^3}$).

Moreover, for $\ell = 13$ these computations and the classification (section 4) show that the image of $\mathbb{P}(\bar{\sigma}_{S_2,13})$ is $\text{PSL}(3, \mathbb{F}_{13})$. The determination of the image for a single prime using the classification appears in [A-M] and [A-A-C].

We have proved the following:

**Theorem 6.1** *For every $\ell \equiv 1 \pmod{4}$, $\ell \not\equiv 1 \pmod{64}$, $\ell \neq 5$, the Galois representation $\bar{\sigma}_{S_2,\lambda}$ is irreducible over $\mathbb{F}_{\ell^3}$. The same holds for $\mathbb{P}(\bar{\sigma}_{S_2,\lambda})$ and $\overline{\chi^{-1} \otimes \sigma_{S_2,\lambda}}$.*

This excludes cases a) and b) of the classification, the case of a diagonal group contained in a maximal subgroup in case c), and the case of a cyclic group contained in a maximal subgroup in case d).

The fact that, after excluding a set of primes $\mathcal{P}$, we obtain just a finite number of reducible primes, is general. The only condition that a compatible family of geometric Galois representations (such as those we obtain from the surfaces $S_a, a \in \mathbb{Q}$) has to verify for this to hold is the following:

**Condition 1.** *If $c$ is (a uniform bound for) the conductor of the family $\sigma_{S_a,\lambda}$, then there exists a prime $p \equiv 1 \pmod{c}$ such that $a_p \notin \mathbb{R}$.*

If we do not know $c$ (or even that it exists) we exclude a set of primes $\mathcal{P}$ such that for every prime $\ell$ decomposing totally in $K := \mathbb{Q}(\{a_p\}) = \mathbb{Q}(i)$ and $\ell \notin \mathcal{P}$ the conductor of $\varepsilon$ (a character ramifying only at the primes of bad reduction of $S_a$ with values in $\mathbb{F}_{\ell^3}^*$) is uniformly bounded by $d$. This is exactly what we have done in the example given by the surface $S_2$ to obtain the "artificial bound" 128 for the conductor. Then we require:

**Condition 1'.** There exists a prime $p \equiv 1 \pmod{d}$ with $a_p \notin \mathbb{R}$.

Under condition 1 or 1', we prove that there is a finite number of reducible primes, outside $\mathcal{P}$ and decomposing totally in $K$, as follows:
Take $p \equiv 1 \pmod{d}$ (or $p \equiv 1 \pmod{c}$) with $a_p \notin \mathbb{R}$.
Then $\varepsilon(p) = 1$ and putting $\mu = \varepsilon \chi^i$, $i = 0, 1, 2$ we obtain:

$$p^3 - 1 \equiv p\bar{a}_p - a_p \pmod{\ell} \quad (i = 0),$$
$$p^3 - 1 \equiv pa_p - \bar{a}_p \pmod{\ell} \quad (i = 2),$$
$$a_p - \bar{a}_p \equiv 0 \pmod{\ell} \quad (i = 1).$$

For the finiteness of the set of solutions of the third, use condition 1 or 1'.
For the other two congruences, use $|a_p| \leq 3p$.

## 6.2 Cases c) and d) of the classification

We will now study the non-cyclic case of d) and case c) of the classification (see section 4).

d) Let $G_\lambda := \mathrm{Image}(\bar{\sigma}_{S_a, \lambda})$. If its projective image $\mathbb{P}(G_\lambda)$ is an irreducible group contained in a maximal subgroup in case d), then $\mathbb{P}(G_\lambda)$ is contained in the normalizer of a cyclic group, and it fits in the exact sequence

$$1 \to \mathcal{C} \to \mathbb{P}(G_\lambda) \to C_3 \to 1,$$

where $\mathcal{C}$ is cyclic of order $|\mathcal{C}| \mid p^2 + p + 1$.

c) If $\mathbb{P}(G_\lambda)$ is an irreducible subgroup contained in a maximal group in case c), then $\mathbb{P}(G_\lambda)$ is contained in the normalizer of a diagonal group, and it fits in the exact sequence

$$1 \to \mathcal{R} \to \mathbb{P}(G_\lambda) \to \mathcal{U} \to 1,$$

where $\mathcal{U}$ is isomorphic to a subgroup of $S_3$ and $\mathcal{R}$ is diagonal of order $|\mathcal{R}| \mid (p-1)^2$.

• **Case d):** Using the description of the image of the inertia subgroup at $\ell$ (see section 5), we know that its projective image will be contained in $\mathcal{C}$, because it is cyclic with order greater or equal to $\ell - 1 > 3$, whenever $\ell \geq 5$.
Then composing $\mathbb{P}(\bar{\sigma}_{S_a, \lambda})$ with the quotient $\mathbb{P}(G_\lambda) \to C_3$ we obtain a representation

$$G_\mathbb{Q} \twoheadrightarrow C_3$$

unramified outside $N$, because when we take quotient by $\mathcal{C}$ we are trivializing the image of the inertia group at $\ell$.

So its kernel gives a cubic character $\psi$ unramified outside $N$, and if $p \nmid \ell N$:

$$\psi(p) \neq 1 \Rightarrow \text{trace}(\mathbb{P}(\bar{\sigma}_{S_a,\lambda})(\text{Frob } p)) = 0 \Leftrightarrow a_p \equiv 0 \pmod{\ell}. \quad (6.1)$$

Here we are using the fact that a matrix obtained by a non trivial cyclic permutation of the columns of a diagonal matrix has trace 0. To control this case we need the following

**Definition 6.2** *let $N'$ be the product of primes of bad reduction of $S_a$ except that we put $9$ in this product if $3$ is a prime of bad reduction. Then we say that the family of representations $\sigma_{S_a,\lambda}$ has cubic Complex Multiplication (or cubic CM) if there exists a character:*

$$\psi : (\mathbb{Z}/N'\mathbb{Z})^* \to C_3$$

*such that for every $p \nmid N'$ with $\psi(p) \neq 1$, $a_p = 0$, i.e,*

$$a_p = \psi(p) a_p \quad \text{for every} \quad p \nmid N'.$$

Then, it is clear that only finitely many primes $\ell$ can fall in case d) if we avoid cubic CM.

In the example $\sigma_{S_2,\lambda}$, $N = 2$ and $3 \nmid (\mathbb{Z}/2^k\mathbb{Z})^*$, so that case d) cannot happen for any $\ell$.

• **Case c):** Again we have $\mathbb{P}(\bar{\sigma}_{S_a,\lambda}|_I)$ contained in the center $\mathcal{R}$, because $\ell - 1 > 3$ for $\ell \geq 5$ and $S_3$ does not contain cyclic subgroups of order greater than 3. Assuming irreducibility of the representation, we must have $\mathcal{U} = S_3$ or $C_3$. If $\mathcal{U} = C_3$, the same analysis as in case d) applies.

If $\mathcal{U} = S_3$, consider the quotient $S_3 \twoheadrightarrow C_2$. Composing the epimorphisms

$$G_\mathbb{Q} \twoheadrightarrow \mathbb{P}(G_\lambda) \twoheadrightarrow S_3 \twoheadrightarrow C_2$$

we obtain a quadratic character $\gamma$ unramified outside $N$. We have that for every $p \nmid \ell N$, $\gamma(p) = -1$ implies that $\mathbb{P}(\bar{\sigma}_{S_a,\lambda})(\text{Frob } p)$ becomes diagonal after an odd permutation of its columns.

These matrices have the following property: If their characteristic polynomial is $x^3 + Ax^2 + Bx + C$, then $AB = C$.

But the characteristic polynomial of $\bar{\sigma}_{S_a,\lambda}(\text{Frob } p)$ is, for every $p \nmid N\ell$:

$$x^3 - a_p x^2 + p\bar{a}_p x - p^3 \pmod{\ell}.$$

Then $\gamma(p) = -1 \Rightarrow -p a_p \bar{a}_p \equiv -p^3 \pmod{\ell}$ must hold, and therefore:

$$a_p \bar{a}_p \equiv p^2 \pmod{\ell}. \quad (6.2)$$

**Definition 6.3** *We say that the family of representations $\sigma_{S_a,\lambda}$ is dual with respect to a quadratic character $\gamma$ (ramifying only at the primes of bad reduction of $S_a$) if, for every prime $p$,*

$$\gamma(p) = -1 \Rightarrow a_p \bar{a}_p = p^2.$$

Again, it is clear that, assuming that $\sigma_{S_a,\lambda}$ is not dual with respect to any quadratic character, it follows that only finitely many primes $\ell$ can fall in case c).

Remark: The consideration of the quotient $C_2$ of $S_3$ is just to simplify computations. Note that the condition obtained is sufficient to guarantee that only finitely many primes fall in case c), but not necessary. The alternative option is to work with $S_3$, and we obtain a situation analogous to that encountered in case d). The difference is that now we have to work with the closure $L$ of a non-cyclic cubic extension of $\mathbb{Q}$, instead of with a cyclic cubic field. $L$ only ramifies at primes of bad reduction of $S_a$ (only finitely many possible $L$ thanks to Hermite's theorem). For every $p \nmid \ell N$ we should look at the element in $\mathrm{Gal}(L/\mathbb{Q})$ corresponding to Frob $p$ : if it is a non trivial cyclic permutation then formula (6.1) applies, and (6.2) if is is an odd permutation. Therefore, for case c) to hold for infinitely many primes the representations should be dual respect to a quadratic character and the coefficients $a_p$ should be 0 for "one third" of the primes.

In the example $\sigma_{S_2,\lambda}$, $\gamma$ would ramify only at 2. Then $\gamma(3) = -1$ or $\gamma(5) = -1$. Using only $a_3 = 1 + 2i$ and $a_5 = -1 - 4i$ we conclude that the representation is not dual respect to $\gamma$. Moreover:
$a_3 \bar{a}_3 = 5 \equiv 3^2 \pmod{\ell} \Rightarrow \ell = 2$,
$a_5 \bar{a}_5 = 17 \equiv 5^2 \pmod{\ell} \Rightarrow \ell = 2$.
So that no prime $\ell \equiv 1 \pmod 4$, $\ell > 5$, falls in case c) .

## 6.3   Case e) : The stabilizer of a conic

Choosing an adequate basis we can write the equation of a general conic as:

$$y^2 = kxz.$$

Then, it is easy to see that the stabilizer of this conic is the group of matrices of the form (cf. [V-Y]):

$$\begin{pmatrix} k\alpha^2 & 2\alpha\beta & \beta^2 \\ k\sqrt{k}\alpha\beta & \sqrt{k}(\alpha\delta + \beta\gamma) & \sqrt{k}\beta\delta \\ k\gamma^2 & 2\gamma\delta & \delta^2 \end{pmatrix}.$$

One can check by direct computation that the characteristic polynomial is for such matrices:

$$x^3 - Ax^2 + Bx - C,$$

with $B^3 = CA^3$.

Therefore, if we suppose that the image of $\bar{\sigma}_{S_a,\lambda}$ falls in case e) we must have: $p^3 \bar{a}_p^3 \equiv p^3 a_p^3 \pmod{\lambda}$ for every $p \nmid \ell N$, or equivalently:

$$a_p^3 \equiv \bar{a}_p^3 \pmod{\lambda}.$$

The values $a_p$ are all in a quadratic imaginary field $K$. Then, for a prime $p$ such that $a_p \notin \mathbb{R}$, $a_p^3 \neq \bar{a}_p^3$; provided that $K \neq \mathbb{Q}(\sqrt{-3})$.

Thus, for the representations $\bar{\sigma}_{S_a,\lambda}$ case e) can only happen for finitely many primes.

More generally, this will be the case whenever the representations are non-selfdual in the more general sense that there does not exist any character $\nu$ such that $a_p = \nu(p)\bar{a}_p$ for every $p \nmid \ell N$.

In the example, using only $a_3$ we obtain: $(1 + 2i)^3 \equiv (1 - 2i)^3 \pmod{\lambda}$, and this implies $\ell = 2$. We conclude that no prime $\ell \equiv 1 \pmod{4}$ falls in case e) in our example.

## 6.4 The exceptional cases

In these cases $\mathbb{P}(G_\lambda) \subseteq \mathrm{PSL}(3, \mathbb{F}_\ell)$ is a group contained in another group of order $36, 72, 168, 216$ or $360$. The description of $\mathbb{P}(\bar{\sigma}_{S_a,\lambda}|_I)$ given in section 5 tells us that in the image of $\mathbb{P}(\bar{\sigma}_{S_a,\lambda})$ there is a cyclic subgroup of order $\ell - 1$, $\ell + 1$ or $\ell^2 + \ell + 1$. This gives only a few possibilities for the exceptional primes. Using the structure of the exceptional groups in cases h), i) and j) and the restriction to $\ell \equiv 1 \pmod{4}$ we end up with $\ell = 5$ or $13$.

As we already mentioned, we know that $\mathbb{P}(\bar{\sigma}_{S_2,13})$ has image $\mathrm{PSL}(3, \mathbb{F}_{13})$. Then, 5 is the only prime $\ell \equiv 1 \pmod{4}$ not in $\mathcal{P}$ such that the image of $\mathbb{P}(\bar{\sigma}_{S_2,\lambda})$ for $\lambda \mid \ell$ can be a proper subgroup of $\mathrm{PSL}(3, \mathbb{F}_\ell)$.

## 6.5 Conclusion

We have considered all cases in the classification of subgroups of $\mathrm{PSL}(3, \mathbb{F}_\ell)$ given in section 4. Therefore, we conclude:

**Theorem 6.4** *The family of representations $\sigma_{S_2,\lambda}$ verifies:*

$$\mathrm{Image}(\mathbb{P}(\bar{\sigma}_{S_2,\lambda})) = \mathrm{PSL}(3, \mathbb{F}_\ell),$$

*for every prime $\ell \equiv 1 \pmod{4}$, $\ell \not\equiv 1 \pmod{64}$, $\ell > 5$.*

Remark: The groups $\mathrm{PSL}(3, \mathbb{F}_\ell)$, for all primes $\ell \equiv 1 \pmod{4}$, have been previously realized as Galois groups over $\mathbb{Q}(t)$ using rigidity methods in [Th].

Observe that the image of $\overline{\chi^{-1} \otimes \sigma_{S_2,\lambda}}$ is contained in $\mathrm{SL}(3, \mathbb{F}_\ell)$ and

$$\mathbb{P}(\overline{\chi^{-1} \otimes \sigma_{S_2,\lambda}}) = \mathbb{P}(\bar{\sigma}_{S_2,\lambda}).$$

Then we conclude that:

**Corollary 6.5**
$$\text{Image}(\overline{\chi^{-1} \otimes \sigma_{S_2,\lambda}}) = \text{SL}(3, \mathbb{F}_\ell),$$

for every $\ell \equiv 1 \pmod 4$, $\ell \not\equiv 1 \pmod{64}$, $\ell > 5$. In particular, for these values of $\ell$ the groups $\text{SL}(3, \mathbb{F}_\ell)$ are Galois groups over $\mathbb{Q}$. For each $\ell$, the corresponding number field is unramified outside $2\ell$.

Remark: Using a lemma of Serre (see [S1]), we also conclude that for these $\ell$ the image of $\chi^{-1} \otimes \sigma_{S_2,\lambda}$ is $\text{SL}(3, \mathbb{Z}_\ell)$.

# 7 On the density of the exceptional set - Condition 1' revisited

We have proved that the image of $\bar{\sigma}_{S_2,\lambda}$ is "as large as possible" for infinitely many totally decomposed primes. Now we will show that this, together with the Cebotarev density theorem, forces the image of the rest of the totally decomposed primes to be also "as large as possible", except at most for a set of Dirichlet density 0.

**Lemma 7.1** *If a compatible family of three-dimensional geometric Galois representations $\sigma_\lambda$ verifies*

$$\text{Image}(\bar{\sigma}_\lambda) = \text{PSL}(3, \mathbb{F}_\ell)$$

*for an infinite set $\mathbb{L}$ of totally decomposed primes, then condition 1' is satisfied for any d.*

Before proving this lemma, let us show its consequences. We denote by $\mathcal{D}$ the set of totally decomposed primes in $\mathbb{Q}(\{a_p\})$. Suppose that for a compatible family of geometric Galois representations we have checked all conditions, as in the previous section, to ensure an image "as large as possible" for those primes in $\mathcal{D}$ that are not in a finite exceptional set $\mathcal{S}$ nor in a certain excluded set $\mathcal{P}$ as in condition 1', and assume that the Dirichlet density of $\mathcal{P}$ is smaller than that of $\mathcal{D}$. Then we apply lemma 7.1 and conclude that condition 1' is still valid for an arbitrarily large $d$ and the corresponding set $\mathcal{P}'$ of arbitrarily small Dirichlet density to be excluded. Therefore, the image must be "as large as possible" for every prime in $\mathcal{D}$ except for the elements of $\mathcal{P}'$ and a finite set $\mathcal{S}'$. Note that $\mathcal{S}'$ may be larger than $\mathcal{S}$, but it is still finite (condition 1' implies that there are only finitely many primes with reducible image in $\mathcal{D} \smallsetminus \mathcal{P}'$). Letting the conductor $d$ go to infinity we conclude that the exceptional set has Dirichlet density 0.

Applying this to the representations $\sigma_{S_2,\lambda}$ we obtain:

**Theorem 7.2**
$$\text{Image}(\overline{\chi^{-1} \otimes \sigma_{S_2,\lambda}}) = \text{SL}(3, \mathbb{F}_\ell)$$

*for every $\ell \equiv 1 \pmod 4$, except at most for a set of Dirichlet density 0.*

*Proof of lemma 7.1.*

The proof is based on the fact that, whenever $a_p$ is real, the characteristic polynomial of $\sigma_\lambda(\text{Frob } p)$ decomposes over $\mathbb{Q}$:

$$x^3 - a_p x^2 + p a_p x - p^3 = (x-p)(x^2 + (p-a_p)x + p^2).$$

Suppose that this happens for every $p \equiv 1 \pmod{d}$ for some $d$. Now consider the residual representation $\overline{\chi^{-1} \otimes \sigma_\lambda}$ for a prime $\ell \in \mathbb{L}$. For such a prime the image is $\text{SL}(3, \mathbb{F}_\ell)$. Let $\mathcal{U}_\ell$ be the image of the subgroup of $\text{Gal}(\bar{\mathbb{Q}}/\mathbb{Q})$ topologically generated by the Frobenius elements at primes $p \equiv 1 \pmod{d}$; it is a subgroup of $\text{SL}(3, \mathbb{F}_\ell)$ and it follows from the Cebotarev density theorem that it is made of (not only generated by) matrices whose characteristic polynomials are reducible over $\mathbb{F}_\ell$. Therefore, it is clear that

$$|\mathcal{U}_\ell|/|\text{SL}(3, \mathbb{F}_\ell)| \leq 1/(\ell^2 + \ell + 1).$$

On the other hand, by the Cebotarev density theorem and the fact that the Dirichlet density of the set of $p \equiv 1 \pmod{d}$ is $1/\varphi(d)$ (where $\varphi$ is the Euler function), we must have:

$$|\mathcal{U}_\ell|/|\text{SL}(3, \mathbb{F}_\ell)| \geq 1/\varphi(d).$$

Therefore $\ell^2 + \ell + 1 \leq \varphi(d)$. But $\ell$ can be taken to be arbitrarily large, because $\mathbb{L}$ is infinite, so we obtain a contradiction, and this proves the lemma.

# 8 Determination of the images: inert primes

## 8.1 Imitation of section 6

With slight changes, the same ideas used for totally decomposed primes can also be applied to deal with the same cases of the classification for inert primes. In our example these are the primes $\ell \equiv 3 \pmod 4$ and for them the image of the representation $\bar{\sigma}_{S_2,\lambda}$ is a subgroup of $\text{PSL}(3, \mathbb{F}_{\ell^2})$. We briefly describe this procedure:

• Reducible representations: We have a character $\mu = \chi^i \varepsilon$ with $\mu(p)$ a root of the characteristic polynomial of $\bar{\sigma}_{S_2,\lambda}(\text{Frob } p)$, $i = 0, 1, 2$ and $\varepsilon$ unramified outside 2, as in section 6.1. The difference is that the image of $\varepsilon$ is now contained in $\mathbb{F}_{\ell^6}^*$. To bound the conductor we subtract from the set of inert primes the set $\mathcal{P}$ of primes $\ell \equiv -1 \pmod{32}$. We obtain $64 \nmid \ell^6 - 1$, and the same bound for the conductor as in section 6.1. So with the computations of that section, using $\varepsilon(127) = 1$ or $\varepsilon(3) = \pm 1$, and the characteristic polynomials corresponding to these two primes, we prove that the image is irreducible (over $\mathbb{F}_{\ell^6}$) for every $\ell \equiv 3 \pmod{4}$, $\ell \not\equiv -1 \pmod{32}$, except of course for $\ell = 3$ and 127. To save the prime 127 we compute the reduction mod 127 of a few characteristic polynomials and this shows that the image of the corresponding representation is irreducible over $\mathbb{F}_{127^2}$ and is not cyclic. We ignore the prime 3 because our method fails for it in cases c) and d), and by direct inspection we cannot determine whether it falls in the exceptional cases or not; so we will assume henceforth that we are dealing only with primes greater than 3.

We have proved that for every $\ell \equiv 3 \pmod 4$, $\ell \not\equiv -1 \pmod{32}$, $\ell > 3$, cases a) , b) and the reducible subcases of c) and d) in the classification are not possible for the image of the corresponding residual Galois representation.

• Cases c), d) and e): All we have done in sections 6.2 and 6.3 applies as well to the inert primes greater than 3 (this inequality being essential in the argument used in cases c) and d)). We conclude that, under the hypothesis of irreducibility, no inert prime greater than 3 falls in these cases (recall that in the computations we have only used $a_3$ and $a_5$).

• Exceptional cases: As in section 6.4, we have a cyclic subgroup of the image with order $\ell - 1$, $\ell + 1$ or $\ell^2 + \ell + 1$. Using the structure of the exceptional groups, we see that they have a cyclic subgroup of one of these orders only for $\ell = 7, 11$ and $19$, among all primes $\ell \equiv 3 \pmod 4$, $\ell > 3$. Computing the reduction of some characteristic polynomials modulo these three primes, we see that there are elements in the image of the corresponding representation whose orders do not divide the order of those exceptional groups we were suspecting the image could be equal to. We conclude that no inert prime greater than 3 falls in the exceptional cases.

Now we consider the remaining cases in the classification (see section 4):
Case f) with $k = 1$: If, for an inert prime $\ell$, the image were contained in $\mathrm{PSL}(3, \mathbb{F}_\ell)$, there should exist for every prime $p \nmid 2\ell$ an element $c(p) \in \bar{\mathbb{F}}_\ell^*$ such that $c(p)^3 p^3 \in \mathbb{F}_\ell^*$ and $c(p) a_p \in \mathbb{F}_\ell$. Therefore, $a_p^3 \in \mathbb{F}_\ell$. This implies that the prime $\ell$ falls in case e), the stabilizer of a conic.

In his article [M], Mitchell also classifies the subgroups of $\mathrm{PSU}(3, \mathbb{F}_\ell)$. All of them are subgroups of those appearing in section 4, so we can conclude from the computations above that no proper subgroup of $\mathrm{PSU}(3, \mathbb{F}_\ell)$ will be obtained as the image of $\sigma_{S_2, \lambda}$ for any inert prime $\ell$ not in the excluded set.

Applying the classification in section 4, for $t = 2$, we obtain:

**Theorem 8.1** *The family of representations $\sigma_{S_2, \lambda}$ verifies:*

$$\mathrm{Image}(\mathbb{P}(\bar{\sigma}_\lambda)) = \mathrm{PSU}(3, \mathbb{F}_\ell) \text{ or } \mathrm{PSL}(3, \mathbb{F}_{\ell^2})$$

*for every prime $\ell \equiv 3 \pmod 4$, $\ell \not\equiv -1 \pmod{32}$, $\ell > 3$.*

Applying lemma 7.1, and letting the bound of the conductor of the character in the reducible case tend to infinity, we conclude, as in section 7, that the density of the set of inert primes with reducible image (over $\mathbb{F}_{\ell^6}$) is 0. For those with irreducible image, the previous theorem applies:

**Corollary 8.2** *The conclusion of theorem 8.1 is valid for every $\ell \equiv 3 \pmod 4$, except at most for a set of Dirichlet density $0$.*

## 8.2 Unitariness

In this section we will prove that for those inert primes such that theorem 8.1 or corollary 8.2 applies, the image is in fact unitary. Our reasoning will be by reductio ad absurdum.

**Lemma 8.3** *For the family of representation $\sigma_{S_a,\lambda}$, the image of $\mathbb{P}(\bar{\sigma}_{S_a,\lambda})$ can not be $\mathrm{PSL}(3, \mathbb{F}_{\ell^2})$ for any prime $\ell \equiv 3 \pmod 4$.*

This lemma, together with theorem 8.1 and corollary 8.2 gives:

**Theorem 8.4** *The family of representations $\sigma_{S_2,\lambda}$ verifies:*

$$\mathrm{Image}(\mathbb{P}(\bar{\sigma}_{S_2,\lambda})) = \mathrm{PSU}(3, \mathbb{F}_\ell)$$

*for every prime $\ell \equiv 3 \pmod 4$, except at most for a set of Dirichlet density $0$, in particular the result is valid for every prime $\ell \equiv 3 \pmod 4$, $\ell \not\equiv -1 \pmod{32}$, $\ell > 3$.*

*Proof of lemma 8.3.*
Suppose that there is a prime $\ell \equiv 3 \pmod 4$ such that the image of $\mathbb{P}(\bar{\sigma}_{S_a,\lambda})$ is $\mathrm{PSL}(3, \mathbb{F}_{\ell^2})$. In particular, there will be elements in this image of order $\ell^4 + \ell^2 + 1$ (projective classes of generators of the cyclic group of units of $\mathbb{F}_{\ell^6}$ embedded in $\mathrm{SL}_3$ via the regular representation). By the Cebotarev density theorem there will be infinitely many primes $p$ such that $\mathbb{P}(\bar{\sigma}_{S_a,\lambda})(\mathrm{Frob}\ p)$ is an element of this order. In this case, the characteristic polynomial of $\bar{\sigma}_{S_a,\lambda}(\mathrm{Frob}\ p)$ will be irreducible modulo $\ell$ and its roots will be elements of $\mathbb{F}_{\ell^6}$ of order multiple of $\ell^4 + \ell^2 + 1$. Therefore the polynomial:

$$x^3 - a_p x^2 + p \bar{a}_p x - p^3 \qquad (8.1)$$

is irreducible in characteristics $0$ and $\ell$, with $a_p \in \mathbb{Q}(i) \smallsetminus \mathbb{Q}$ and its reduction in $\mathbb{F}_{\ell^2} \smallsetminus \mathbb{F}_\ell$.
Let $\alpha$ be one of its roots, it is an algebraic integer with $[\mathbb{Q}(\alpha) : \mathbb{Q}] = 6$ and $\ell$ is inert in $\mathbb{Q}(\alpha)/\mathbb{Q}$. Let $\mathcal{L}$ be the prime in $\mathbb{Q}(\alpha)$ dividing $\ell$ and let us denote by ˆ reduction mod $\mathcal{L}$. Then $[\mathbb{F}_\ell(\hat{\alpha}) : \mathbb{F}_\ell] = 6$.
By reducing mod $\mathcal{L}$, we can see that

$$\alpha^{\ell^3+1} \equiv p^2 \pmod \ell.$$

So the order of $\alpha$ divides $(\ell^3 + 1)(\ell - 1) = (\ell^2 - \ell + 1)(\ell + 1)(\ell - 1)$.
But $\alpha$ is a root of the characteristic polynomial (8.1) and, by assumption, its order is a multiple of $\ell^4 + \ell^2 + 1 = (\ell^2 + \ell + 1)(\ell^2 - \ell + 1)$. This is a contradiction. This proves that for every prime inert in $[\mathbb{Q}(i) : \mathbb{Q}]$, the image of the corresponding Galois representation cannot be $\mathrm{PSL}(3, \mathbb{F}_{\ell^2})$.

Remark: The fact that the images of the Galois representations for inert primes, when irreducible, are unitary, can be proved directly using the particular form of the characteristic polynomials, which in turn is a consequence of the Riemann hypothesis (see [A-M]).

**Corollary 8.5**
$$\text{Image}(\overline{\chi^{-1} \otimes \sigma_{S_2, \lambda}}) = \text{SU}(3, \mathbb{F}_\ell),$$

for every $\ell \equiv 3 \pmod 4$, $\ell \not\equiv -1 \pmod{32}$, $\ell > 3$. In particular, for these values of $\ell$ the groups $\text{SU}(3, \mathbb{F}_\ell)$ are Galois groups over $\mathbb{Q}$. For each $\ell$, the corresponding number field is unramified outside $2\ell$.

Remark: The groups $\text{SU}(3, \mathbb{F}_\ell)$ for all primes $\ell \equiv 3 \pmod 4$, $\ell > 3$, have been realized as Galois groups over $\mathbb{Q}$ using rigidity methods in [Ma].

## 9 Final comments on the example

By Deligne's corollary to the results of de Jong (cf. [Be], Proposition 6.3.2), we know that in our example an uniform bound for the conductor exists. As a consequence of the results in section 7, condition 1 (see section 6) is satisfied. Therefore the image is as large as possible for almost every prime, i.e., theorem 7.2 and corollary 8.5 hold for almost every prime.

For the prime $\ell = 5$ we have computed the reduction of the characteristic polynomials of $\sigma_{S_2,5}(\text{Frob } p)$ modulo 5 for $2 < p \leq 97$, $p \neq 5$, and we have found the trivial character 1 as a root of them; so we suspect that the residual representation for 5 is reducible with $\varepsilon = 1$ and "i=0". If this were so, after semi-simplification we also obtain a two-dimensional subrepresentation $\pi$ over $\mathbb{F}_5$, unramified outside $2 \cdot 5$, with determinant $\chi^3$ ($\chi$ the mod 5 cyclotomic character). The semi-simplification of the restriction of $\pi$ to the inertia group at 5 (see section 5) is the sum of the characters $\chi$ and $\chi^2$ or $\psi_2^{2+5}$ and $\psi_2^{1+2 \cdot 5}$. The computations also show that $\pi$ is irreducible. If $a_p$ is the trace of $\bar{\sigma}_{S_2,5}(\text{Frob } p)$ and $b_p$ that of $\pi(\text{Frob } p)$ then we have:

$$a_p \equiv 1 + b_p \pmod 5, \qquad (10.1)$$

for every $p \nmid 10$. Applying Serre's conjecture (see [S3]) to $\pi$ we deduce that there should exist a classical modular form $f$ of level a power of 2 and weight 8 whose Fourier coefficients reduced modulo 5 agree with the $b_p$. Looking at the space of newforms $S_8(128)$ we have found a newform verifying this for every $p \leq 97$, thus providing more evidence for the fact that $\bar{\sigma}_{S_2,5}$ is reducible. Observe that the congruence (10.1) relates the coefficients of $f$, an eigenform in $S_8(128)$, with those of an eigenform in $H^3(\Gamma_0(128), \mathbb{C})$.

The values $a_p$ are computed in [vG-T1] in terms of the number of points of $S_2$ and other related surfaces over finite fields. The reducibility at 5 would imply congruences modulo 5 involving these numbers of points, depending only on the value modulo 40 of the characteristic $p$ of the finite field.

For an inert prime $\ell$ greater than 3 such that the image of $\bar{\sigma}_{S_2,\lambda}$ is not as in theorem 8.4, if any, this representation must be reducible and with "i=1". After semi-simplification we would obtain a two dimensional subrepresentation $\pi$ over $\mathbb{F}_{\ell^2}$. It can be proved that Image($\mathbb{P}(\pi)$) must be conjugated in $\mathrm{PGL}(2, \mathbb{F}_{\ell^4})$ to $\mathrm{PSL}(2, \mathbb{F}_\ell)$. In this case, after conjugation and twisting by an appropriate character, the image of $\pi$ would be $\mathrm{SL}(2, \mathbb{F}_\ell)$. To see this we apply the classification of subgroups of $\mathrm{PSL}(2, \mathbb{F}_{\ell^r})$ as in [R1] and [R2] and the techniques of this article; unitariness and the Cebotarev density theorem to prove that the reducibility assumption implies that the three roots of the characteristic polynomials of $\bar{\sigma}_{S_2,\lambda}(\mathrm{Frob}\ p)$ must be in $\mathbb{F}_{\ell^2}$, the techniques of sections 6 and 8 to eliminate the dihedral and special cases for the image of $\mathbb{P}(\pi)$ and the description of the restriction to inertia at $\ell$ to eliminate the possibility that $\pi$ were reducible.

## 10 On the generality of the previous results

We stress that conclusions such as those of sections 6, 7 and 8 asserting that the images are generically "as large as possible" will still be valid if we perform the computations starting from any other similar geometric three-dimensional compatible family of Galois representations defined over an imaginary quadratic field $K$ (for example, those obtained from the surfaces $S_a$ for any $a \in \mathbb{Z} - \{0\}$) as long as the Hodge structure is as in the example, condition 1' is verified with an "excluded set" of primes of relatively small density, the family of representations is non-selfdual in the general sense (this is automatic if $K \neq \mathbb{Q}(\sqrt{-3})$) and does not have cubic CM, nor is dual with respect to any character.

In particular, Jasper Scholten has communicated to us that with the techniques of [Sc] using elliptic surfaces he can produce examples of non-selfdual compatible families of three-dimensional Galois representations defined over $\mathbb{Q}(\sqrt{-3})$. Unfortunately, these examples have not been published and are unavailable. If any of these examples verifies all of the above conditions, then as a consequence we would realize the groups $\mathrm{SL}(3, \mathbb{F}_\ell)$ as Galois groups over $\mathbb{Q}$ for every $\ell \equiv \pm 1 \pmod{12}$, except for a subset of primes of Dirichlet density 0, and the groups $\mathrm{SU}(3, \mathbb{F}_\ell)$ for every $\ell \not\equiv \pm 1 \pmod{12}$ except for a subset of Dirichlet density 0.

The only point where the argument seems to depend on the fact that $K = \mathbb{Q}(i)$ is section 8.2, but it is straightforward to verify that it applies to any imaginary quadratic field. The only step that needs a closer look is the congruence $\hat{i}^{\ell^3} = -\hat{i}$ (in fact an equality even in characteristic 0), a crucial point because it shows that complex conjugation agrees with the Frobenius automorphism. In the general case it becomes $\hat{\xi}^{\ell^3} = -\hat{\xi}$ with $\xi = \sqrt{-D}$ and $\ell$ inert in $K = \mathbb{Q}(\xi)$. This congruence follows from $(-D)^{\frac{\ell-1}{2}} \equiv -1 \pmod{\ell}$.

# 11 Strong Form of Clozel's conjecture and applications

## 11.1 The Conjecture

In this section we will assume the validity of a strong form of Clozel's conjecture and we will derive consequences for inverse Galois theory. We will give intrinsic conditions on a Hecke eigenform (whose eigenvalues generate an imaginary quadratic field) guaranteeing that the results summarized in the previous section apply to the attached Galois representations, so for those eigenforms verifying these conditions the images of the residual representations will be "as large as possible" for almost every prime: $SL_3$ or $SU_3$ (after twisting by $\chi$) depending on the decomposition type of the prime.

**Conjecture 11.1** *Let $f \in H^3(\Gamma_0(N), \mathbb{C})$ be a cuspidal Hecke eigenform, then there is a compatible family of Galois representations attached to $f$ as described in conjecture 1.4. These representations $\rho_\lambda$ arise from the action of $G_\mathbb{Q}$ on a rank 3 motive defined over $\mathbb{Q}$ with coefficients in $\mathbb{Q}_f$. For every $\ell \nmid N$ the (semi-simplification of the) restriction of $\bar\rho_\lambda$ to the inertia group at $\ell$ is not the group of scalar matrices $\chi \cdot Id$. Moreover, for every $\ell$, the conductor of $\bar\rho_\lambda$ divides the level $N$.*

Remark: The resulting motive has weight 2.
Justification: The geometric or motivic nature of the conjectured Galois representations appears already in [C] (see also [Ra]). The relation between the weight of the modular form and the image of inertia is taken from the generalizations of Serre's conjectures by Ash and Sinnott (see [A-Si]).
Finally, the relation between the level and the conductor comes from the analogy with the case of classical modular forms. Recently, this relation has been postulated as part of the generalized Serre's conjecture for the case of arbitrary level (see [A-D-P]). Removing this last hypothesis would complicate the computations as in the purely geometric example we considered above, and the results on the largeness of the images would hold only outside a density 0 set of primes. Assuming conjecture 11.1, we can apply the techniques of the previous sections to modular families of Galois representation. The desired result will follow as long as the modular form $f$ verify certain conditions.

## 11.2 Conditions on the Modular Forms

As usual, we consider only the case of Hecke eigenforms

$$f \in H^3(\Gamma_0(N), \mathbb{C})$$

with $\mathbb{Q}_f$ imaginary quadratic. If $f$ is cuspidal, the needed facts, i.e., that the roots of the characteristic polynomial of the image of Frob $p$ should have absolute value $p$ (Ramanujan's conjecture) and the bound on the eigenvalue $a_p$ that this implies, are known (using the fact that $f$ is a cuspidal automorphic forms on GL(3)). They agree with the geometric hypothesis in conjecture 11.1 (Riemann hypothesis).

The statement about the image of inertia at $\ell$ in conjecture 11.1 can be replaced by a condition on the Hodge decomposition of the motive, as explained in section 5.

Therefore, in order to apply the techniques of previous sections to these representations and obtain large images, we only need to determine which conditions on $f$ correspond to the fact that the family of Galois representations become:

1) selfdual (in general sense),
2) reducible, or
3) monomial.

Explanation:
2) If the representations $\bar{\rho}_\lambda$ were reducible for infinitely many $\lambda$, then thanks to Ramanujan's conjecture we can assume (see sections 6 and 8) that one of the irreducible components is of the form $\varepsilon_\lambda \chi$, where $\varepsilon_\lambda$ is a character with conductor dividing $N$. Moreover, using the fact that the determinant $\chi^3$ of $\bar{\rho}_\lambda$ is unramified at $N$, we conclude that it must hold: $\text{cond}(\varepsilon_\lambda)^2 \mid N$.

Let us show that in this situation the representations are reducible also in characteristic 0, independently of $\ell$. To see this, observe that from the uniform bound on the conductors of the characters $\varepsilon_\lambda$, we know that there is a complex character $\varepsilon$ with the square of its conductor dividing $N$ such that $\varepsilon(p)p$ is a root of $\rho_\lambda(\text{Frob } p)$, for every $p \nmid N$, $\ell \nmid N$.

From this and the reducibility of infinitely many residual representations, an application of the theory of pseudo-representations (see [T]) proves that the values $\{a_p - \varepsilon(p)p\}$ correspond to the traces of a two-dimensional compatible family of Galois representations $\{\nu'_\lambda\}$.

For an exposition of the standard facts of the theory of pseudo-representations and how to use them to prove reducibility, see [D1], chapter 6 (or [D2]). The only difference is that because we are not dealing with odd residual representations the construction of the representations $\nu'_\lambda$ is done using Taylor's notion of pseudo-representation instead of Wiles' one. Therefore, for every $\lambda$ there is a $\lambda'$-adic Galois representation $\nu_{\lambda'}$ over $\overline{\mathbb{Q}}_{\lambda'}$ ($\lambda'$ a prime above $\lambda$) with $\text{trace}(\nu_{\lambda'}(\text{Frob } p)) = a_p - \varepsilon(p)p$, for every $p \nmid \ell N$.

We conclude that if the representations $\bar{\rho}_\lambda$ are reducible for infinitely many $\lambda$, then they are so for every prime $\lambda$. Furthermore, after a field extension we obtain a reducible family of $\lambda'$-adic representations:

$$\rho^{ss}_{\lambda'} = \varepsilon \cdot \chi \oplus \nu_{\lambda'}.$$

The converse is obvious. In what follows, condition 2) refers to either of these two equivalent conditions.

3) As explained in section 6.2, dual respect to a character was a weaker notion introduced to simplify computations, but if the images of the Galois representations fall in case c) or d) (see section 4) for infinitely many primes, the representations are monomial for every $\lambda \nmid N$. In these cases the representations are induced from a character of a cyclic, respectively (Galois closure of a) non-cyclic, cubic extension of $\mathbb{Q}$.

1) "In general sense" means up to multiplication by a character. Recall that this is incompatible with the assumption that $\mathbb{Q}_f$ is imaginary quadratic, except for the case of $\mathbb{Q}_f = \mathbb{Q}(\sqrt{-3})$.

Now let us translate these three conditions to conditions on the Hecke eigenform $f$, using the language of automorphic forms. This can be done because our space of modular forms can be thought of as being included in that of automorphic forms on $GL(3)/\mathbb{Q}$, and the cuspidal modular forms correspond to cuspidal automorphic forms (see [A-G-G]).

**Condition 1)** implies that the automorphic form is selfdual:

$$f \otimes \omega = f^{\vee}, \qquad (I)$$

where $\omega$ is a Hecke character, i.e., a one dimensional automorphic representation of $GL(1)/\mathbb{Q}$. This is the case in particular when $f$ is the symmetric square constructed by Gelbart and Jacquet of an automorphic representation on $GL(2)/\mathbb{Q}$.

**Condition 2)** implies that the Galois representations attached to $f$ (assuming they are semisimple) reduce as the sum of the Galois representations associated to a Hecke character and those associated to a certain rank two motive. In this situation, we can apply a result ([Ra], prop. 3.2.4) guaranteeing that the property of a motive of being associated to an automorphic form is "closed under subtraction" (this generalizes theorem 4.7 in [J-S]), so we have:

$$f = \omega \boxplus g, \qquad (II)$$

where $\omega$ is a Hecke character and $g$ an automorphic form on $GL(2)/\mathbb{Q}$. In this case $f$ cannot be cuspidal (see [J-S]).

Remark: The result we are applying is only proved (see [Ra]) under the assumption of the analytic Tate conjecture.

**Condition 3)** in our case of cubic induction is also a case of Langlands' functoriality solved, i.e., given $E$ (the Galois closure of) a cubic extension of $\mathbb{Q}$ and a Hecke character $\varphi$ on $GL(1)/E$ there is an automorphic form on $GL(3)/\mathbb{Q}$ with the same $L$-function:

$$f = Ind_E^{\mathbb{Q}}(\varphi). \qquad (III)$$

These are the automorphic representations corresponding to condition 3).

**Theorem 11.2** *Assume conjecture 11.1 and the analytic Tate conjecture (cf. [Ra]). If a cuspidal Hecke eigenform $f$ with $\mathbb{Q}_f$ imaginary quadratic does not verify either condition (I) or (III) above, then the images of the residual representations $\bar{\rho}_\lambda$ attached to $f$ are "as large as possible" for almost every $\ell$. Moreover, if $\mathbb{Q}_f \neq \mathbb{Q}(\sqrt{-3})$ then $f$ does not verify condition (I).*

## 11.3 Examples

First we apply these ideas to the eigenform $f \in H^3(\Gamma_0(128), \mathbb{C})$ appearing in theorem 2.1. As stated there, all coefficients used in the computations agree with the traces of the family of geometric Galois representations studied in detail in previous sections. In particular from the computations done in those sections we conclude that this eigenform does not verify condition (I), (II) or (III). The artificial conductor for condition 1') used in section 6 and 8 was 128, which we are now assuming to be the conductor of the family of Galois representations attached to $f$, so conditions 1') and 1) become equivalent in this case. Therefore, the conclusions of sections 6 and 8 follow, but now for all primes.

**Theorem 11.3** *Assume conjecture 11.1 to be true. Then, if we call $\sigma_\lambda$ the family of Galois representations attached to the eigenform $f \in H^3(\Gamma_0(128), \mathbb{C})$, we have*
$$\mathrm{Image}(\overline{\chi^{-1} \otimes \sigma_\lambda}) = \mathrm{SL}(3, \mathbb{F}_\ell),$$
*for every $\ell \equiv 1 \pmod 4$, $\ell > 5$ and*
$$\mathrm{Image}(\overline{\chi^{-1} \otimes \sigma_\lambda}) = \mathrm{SU}(3, \mathbb{F}_\ell)$$
*for every $\ell \equiv 3 \pmod 4$, $\ell > 3$. In particular, these linear and unitary groups are Galois groups over $\mathbb{Q}$.*

As our next example we take the eigenform $f \in H^3(\Gamma_0(88), \mathbb{C})$ with $\mathbb{Q}_f = \mathbb{Q}(\sqrt{-7})$ whose coefficients are listed in [vG-K-T-V] for $p \leq 173$.
We assume that the conductor of the family of Galois representations attached to $f$ is 88, so condition 1) defined in section 6 is satisfied because $a_{89} = -60 - 4\sqrt{-7}$ is not real. This implies in particular that $f$ does not verify condition (II). Using $a_{89}$ we compute the finitely many primes with (possibly) reducible residual image as described in sections 6 and 8 and we only obtain $\ell = 2, 2879, 48889$ (except for 7 which ramifies in $\mathbb{Q}_f$, 11 which may divide the conductor of the representations and 89). Looking at the reduction of a couple of characteristic polynomials modulo $\ell = 89, 2879, 48889$ we see that for none of them the residual representation is reducible.

Remark: In order to simplify the determination of the primes with reducible image, condition 1) can be improved by using the fact that $\mathrm{cond}(\varepsilon)^2 \mid N$. In particular if $N$ is square-free the character $\varepsilon$ becomes trivial. We will use this in the next example.
As explained in section 6, we study cases c) and d) considering all possible cubic (none!) and quadratic characters unramified outside $2 \cdot 11$. Using a couple of coefficients of $f$ for each possible such character and using the equations described in section 6, we see that no prime $\ell > 3$ falls in these cases, i.e., the image of the residual representation is not monomial for any such $\ell$. In particular, this implies that $f$ is not an eigenform verifying condition (III).
Using only $a_3 = -1 + \sqrt{-7}$ we study condition e) as in section 6 and see that no $\ell > 3$ falls in this case, i.e., the image of the residual representation is not selfdual (in general sense) for any such $\ell$. It was a priori known that $f$ is not an eigenform verifying condition (I) from the last assertion in theorem 11.2.

Finally, as in section 6 and 8, we know from the possible orders of the image of inertia (thanks to the exclusion of the trivial case in conjecture 11.1) that the only primes $\ell > 3$ such that the image may fall in one of the exceptional cases in the classification are: $\ell = 5, 7, 11, 13, 19$. The prime 7 ramifies in $\mathbb{Q}_f$ and 11 may divide the conductor. For the remaining three primes, we compute the reduction of some characteristic polynomials and see that there are elements in the images with orders not corresponding to the order of any element in any of the exceptional groups. So we conclude that the image of the residual Galois representation is not exceptional for any $\ell > 3$.

We remark that for the prime $\ell = 5$, inert in $\mathbb{Q}(\sqrt{-7})$, we have considered (and eliminated) the cases of exceptional image $2 \cdot A_6$ or $A_7$ that were not listed in section 4 and are the only extra possibilities (see [M]).

From all this we conclude:

**Theorem 11.4** *Assume conjecture 11.1 to be true. Then, if we call $\sigma_\lambda$ the family of Galois representations attached to $f$, we have*

$$\mathrm{Image}(\overline{\chi^{-1} \otimes \sigma_\lambda}) = \mathrm{SL}(3, \mathbb{F}_\ell)$$

*for every $\ell$ such that $(\frac{-7}{\ell}) = 1$, $\ell \geq 13$ and*

$$\mathrm{Image}(\overline{\chi^{-1} \otimes \sigma_\lambda}) = \mathrm{SU}(3, \mathbb{F}_\ell)$$

*for every $\ell$ such that $(\frac{-7}{\ell}) = -1$, $\ell \geq 5$.*

*In particular, these linear and unitary groups are Galois groups over $\mathbb{Q}$.*

Combining the last two theorems we also have:

**Corollary 11.5** *Assume conjecture 11.1 to be true. Then, the following groups are Galois groups over $\mathbb{Q}$:*

$$\mathrm{SL}(3, \mathbb{F}_\ell)$$

*for every $\ell \equiv 1, 2, 4 \pmod 7$, $\ell \geq 13$; and*

$$\mathrm{SU}(3, \mathbb{F}_\ell)$$

*for every $\ell \equiv 3, 5, 6 \pmod 7$, $\ell \geq 5$.*

Another example: Consider the eigenform $f \in H^3(\Gamma_0(53), \mathbb{C})$ with $\mathbb{Q}_f = \mathbb{Q}(\sqrt{-11})$. To check reducibility we can take the character $\varepsilon$ as trivial, because the level is prime. Therefore, the congruences determining reducible primes introduced in section 6.1 now hold with any $a_p$. This and all other cases of the classification are dealt with as in the previous examples, using just the first few characteristic polynomials. After the computations we find 5 as the only prime greater than 3 (excluding 11 and 53) such that the image seems to be not "as large as possible". For this prime, as in the geometric example, the representation seems to be reducible with the trivial character as one-dimensional constituent (5 decomposes in $\mathbb{Q}_f$).

This example, combined with the one of level 128, gives the corollary:

**Corollary 11.6** *Assume conjecture 11.1. Then the following groups are Galois groups over $\mathbb{Q}$:*

$$\mathrm{SL}(3, \mathbb{F}_\ell)$$

*for every $\ell \equiv 1, 3, 4, 5, 9 \pmod{11}$, $\ell > 5$; and*

$$\mathrm{SU}(3, \mathbb{F}_\ell)$$

*for every $\ell \equiv 2, 6, 7, 8, 10 \pmod{11}$, $\ell \geq 5$.*

# References

[A-A-C] Allison, G., Ash, A., Conrad, E., Galois Representations, Hecke Operators, and the mod-$p$ Cohomology of $\mathrm{GL}(3, \mathbb{Z})$ with Twisted Coefficients, Experimental Math. **7**:4, (1998) 361-390

[A-D-P] Ash, A., Doud, D., Pollack, D., Galois representations with conjectural connections with arithmetic cohomology, (2001) preprint

[A-G-G] Ash, A., Grayson, D., Green, P., Computations of Cuspidal Cohomology of Congruence Subgroups of $\mathrm{SL}(3, \mathbb{Z})$, J. Number Theory **19**, (1984) 412-436

[A-M] Ash, A., McConnell, M., Experimental indications of three- dimensional Galois representations from the cohomology of $\mathrm{SL}(3, \mathbb{Z})$, Experimental Math. **1**:3, (1992) 209-223

[A-Si] Ash, A., Sinnott, W., An analogue of Serre's conjecture for Galois representations and Hecke eigenclasses in the mod-$p$ cohomology of $\mathrm{GL}(n, \mathbb{Z})$, Duke Math. J. **105**, (2000) 1-24

[A-St] Ash, A., Stevens, G., Cohomology of arithmetic groups and congruences between systems of Hecke eigenvalues, J. Reine Angew. Math. **365**, (1986) 192-220

[Be] Berthelot, P., Altrations de varits algbriques [d'aprs A. J. de Jong], in "Sminaire Bourbaki", vol. 1995/96, exp. 815, 273-311; Astrisque 1997

[B] Bloom, D., The subgroups of $\mathrm{PSL}(3, q)$ for odd $q$, Trans. Amer. Math. Soc. **127**, (1967) 150-178

[C] Clozel, L., Motifs et formes automorphes: applications du principe de fonctorialit, in "Automorphic Forms, Shimura Varieties and $L$-functions", Clozel, L., Milne, J. (eds); Proceedings of the Ann Arbor Conference, 77-159, Academic Press 1990

[D1] Dieulefait, L., Modular Galois Realizations of Linear Groups, thesis, Universitat de Barcelona (2001)

[D2] ____ , On the images of the Galois representations attached to genus 2 Siegel modular forms, preprint, 2001

[F-L] Fontaine, J.M., Laffaille, G., Construction de reprsentations $p$-adiques, Ann. Scient. c. Norm. Sup., $4^e$ srie, t. **15**, (1982) 547-608

[F-M] Fontaine, J.M, Messing, W., $p$-adic periods and $p$-adic tale cohomology, in "Currents Trends in Arithmetical Algebraic Geometry (Arcata, Calif., 1985)"; Contemporary Mathematics **67**, (1987) 179-207

[vG-K-T-V] van Geemen, B., van der Kallen, W., Top, J., Verberkmoes, A., Hecke Eigenforms in the Cohomology of Congruences Subgroups of $SL(3, \mathbb{Z})$, Experiment. Math. **6**:2, (1997) 163-174

[vG-T1] van Geemen, B., Top, J., A non-selfdual automorphic representation of $GL_3$ and a Galois representation, Invent. Math. **117**, (1994) 391-401

[vG-T2] ____ , Selfdual and non-selfdual 3-dimensional Galois representations, Compositio Math. **97** (1995) 51–70

[J-S] Jacquet, H., Shalika, J., Euler products and the classification of automorphic forms II, Amer. J. Math. **103**, (1981) 777-815

[Ma] Malle, G., Some unitary groups as Galois groups over $\mathbb{Q}$, J. Algebra **131**, (1990) 476-482

[M] Mitchell, H., Determination of the ordinary and modular ternary linear groups, Trans. Amer. Math. Soc. **12**, (1911) 207-242

[Ra] Ramakrishnan, D., Pure Motives and Automorphic Forms, Proc. Sympos. Pure Math. **55** Part 2, (1994) 411-445

[R1] Ribet, K.A., On $\ell$-adic representations attached to modular forms, Invent. Math. **28**, (1975) 245-275

[R2] ____ , Images of semistable Galois representations, Pacific J. of Math. **181**, (1997)

[Sc] Scholten, J., Mordell-Weil groups of elliptic surfaces and Galois representations, thesis, Rijksuniversiteit Groningen, (2000)

[S1] Serre, J-P., *Abelian $\ell$-adic representations and elliptic curves*, Benjamin 1968

[S2] ____ , Proprits galoisiennes des points d'ordre fini des courbes elliptiques, Invent. Math. **15**, (1972) 259-331

[S3] ____ , Sur les reprsentations modulaires de degr 2 de $Gal(\bar{\mathbb{Q}}/\mathbb{Q})$, Duke Math. J. **54**, (1987) 179-230

[T] Taylor, R., Galois representations associated to Siegel modular forms of low weight, Duke Math. J. **63**, (1991) 281-332

[Th] Thompson, J., $PSL_3$ as Galois groups over $\mathbb{Q}$, in "Proceedings of the Rutgers group theory year, 1983-1984", 309-319, Cambridge U.P. 1984

[V-Y] Veblen, O., Young, J., *Projective geometry*, Blaisdell Publishing Company 1965